\documentclass[11pt,a4paper]{article}
\usepackage[cp1251]{inputenc}
\usepackage{amsfonts, amssymb, amsthm, amsmath}
\usepackage{mathtext}
\usepackage{longtable}
\usepackage[russian]{babel}
\usepackage{color}

\usepackage[left=1cm,right=1cm,top=2cm,bottom=2cm]{geometry}
\usepackage{graphicx}

\renewcommand{\r}{\mathbb R}

\renewcommand{\le}{\leqslant}
\renewcommand{\ge}{\geqslant}

\renewcommand{\phi}{\varphi}

\newcommand{\eqd}{\stackrel{d}{=}}

\newcommand{\pto}{\stackrel{P}{\longrightarrow}}

\title{A note on mixture representations for the Linnik and Mittag-Leffler
distributions and their applications\thanks{Research supported by
the Russian Science Foundation (project 14-11-00364).}}
\author{V.\,Yu.~Korolev\thanks{Faculty of Computational Mathematics and Cybernetics,
   Moscow State University; Federal Research Center ``Informatics and Control, Russian Academy of Sciences; e-mail: vkorolev@cs.msu.ru},
   A.\,I.~Zeifman\thanks{Vologda State University; Federal Research Center ``Informatics and Control, Russian Academy of Sciences;
   e-mail: a\_zeifman@mail.ru}}
\date{}

\begin{document}

\maketitle

{\small

{\bf Abstarct:} We present some product representations for random
variables with the Linnik, Mittag-Leffler and Weibull distributions
and establish the relationship between the mixing distributions in
these representations. The main result is the representation of the
Linnik distribution as a normal scale mixture with the
Mittag-Leffler mixing distribution. As a corollary, we obtain the
known representation of the Linnik distribution as a scale mixture
of Laplace distributions. Another corollary of the main
representation is the theorem establishing that the distributions of
random sums of independent identically distributed random variables
with finite variances converge to the Linnik distribution under an
appropriate normalization if and only if the distribution of the
random number of summands under the same normalization converges to
the Mittag-Leffler distribution.

\section{Introduction}

Usually the Mittag-Leffler and Linnik distributions are mentioned in
the literature together as examples of geometric stable
distributions. Since these distributions are very often pointed at
as weak limits for geometric random sums, there might have emerged a
prejudice that the scheme of {\it geometric} summation is the only
asymptotic setting within which these distributions can be limiting
for sums of independent and identically distributed random
variables. This prejudice is accompanied by the suspicion that
non-trivial ($\delta<1$, $\alpha<2$) Mittag-Leffler and Linnik laws
can be limiting only for sums in which the summands have infinite
variances.

Another obvious reason for which the Mittag-Leffler and Linnik
distributions are often brought together in the literature is the
formal similarity of the Laplace transform of the former and the
Fourier--Stieltjes transform of the latter.

The main aim of the present paper is to study the analytic and
asymptotic relations between these two laws. We will show that
actually the link between these two laws is much more interesting.
Namely, it turns out that the Linnik distribution with parameter
$\alpha$ is a scale mixture of the normal distributions with the
mixing distribution being the Mittag-Leffler law with parameter
$\delta=\alpha/2$. As a corollary of this result we obtain a theorem
establishing that the distributions of random sums of independent
identically distributed random variables {\it with finite variances}
converge to the Linnik distribution under an appropriate
normalization if and only if the distribution of the random number
of summands under the same normalization converges to the
Mittag-Leffler distribution.

Our main tools are mixture representations for the Linnik,
Mittag-Leffler and Weibull distributions. Mixture representations
for the Linnik and Mittag-Leffler laws were the objects of
investigation in \cite{Devroye1990, ErdoganOstrovskii1997,
ErdoganOstrovskii1998, KotzOstrovskii1996, Pakes1992,
Kozubowski1998, Kozubowski1999}. Some of the results of these papers
will be used in what follows and will be presented as lemmas below.
As well, the proofs of our results are based on some new mixture
representations for the Weibull distributions.

However, we are unaware of the results concerning the possibility of
representation of the Linnik distribution as a scale mixture of
normals. Perhaps, the paper \cite{KotzOstrovskii1996} is the closest
to this conclusion and exposes the representability of the Linnik
law as a scale mixture of Laplace distributions with the mixing
distribution written out explicitly.

We develop the results of \cite{KotzOstrovskii1996} and obtain our
main result which is the representation of the Linnik distribution
as a normal scale mixture with the Mittag-Leffler mixing
distribution, thus finding a tight and clear analytical link between
the Linnik and Mittag-Leffler distributions. As a corollary, we have
the representation of the Linnik distribution as a scale mixture of
Laplace distributions obtained in \cite{KotzOstrovskii1996}. Our
proof of this result together with the explicit formula for the
mixing density obtained in \cite{KotzOstrovskii1996} lead to a
by-product corollary of this representation which is the explicit
formula for the distribution density of the ratio of two independent
positive strictly stable random variables. Another consequence of
the main representation is the theorem establishing that the
distributions of random sums of independent identically distributed
random variables with finite variances converge to the Linnik
distribution under an appropriate normalization if and only if the
distribution of the random number of summands under the same
normalization converges to the Mittag-Leffler distribution. On the
one hand, this theorem offers a new <<asymptotic>> link between the
Linnik and Mittag-Leffler distributions and, on the other hand, it
shows that the Linnik law can be limiting for the distribution of
sums of independent identically distributed random variables {\it
with finite variances} thus dispelling the suspicion mentioned
above.

The paper is organized as follows. Section 2 presents the
definitions and basic properties of the Linnik and Mittag-Leffler
distributions. Section 3 contains all the definitions and auxiliary
results. The proofs of our main results are purposely indirect and
essentially rely on some new mixture properties of the Weibull
distribution also presented in Section 3. In Section 4 we prove the
representability of the Linnik distribution as the scale mixture of
normal laws with the Mittag-Leffler mixing distribution and as the
scale mixture of the Laplace laws with the mixing distribution being
that of the ratio of two independent random variables with the same
strictly stable distribution concentrated on the nonnegative
halfline. We use this representation together with the result of
\cite{KotzOstrovskii1996} to obtain a by-product corollary which is
the explicit representation of the distribution density of the ratio
of two independent positive strictly stable random variables. In
Section 5 we prove and discuss the theorem establishing that the
distributions of random sums of independent identically distributed
random variables with finite variances converge to the Linnik
distribution under an appropriate normalization if and only if the
distribution of the random number of summands under the same
normalization converges to the Mittag-Leffler distribution.

\section{The Mittag-Leffler and Linnik distributions}

\subsection{The Mittag-Leffler distributions}

The Mittag-Leffler probability distribution is the distribution of a
nonnegative random variable $M_{\delta}$ whose Laplace transform is
$$
\psi_{\delta}(s)\equiv {\sf E}e^{-sM_{\delta}}=\frac{1}{1+\lambda
s^{\delta}},\ \ \ s\ge0,\eqno(1)
$$
where $\lambda>0$, $0<\delta\le1$. For simplicity, in what follows
we will consider the standard scale case and assume that
$\lambda=1$.

The origin of the term {\it Mittag-Leffler distribution} is due to
that the probability density corresponding to Laplace transform (1)
has the form
$$
f_{\delta}^{ML}(x)=\frac{1}{x^{1-\delta}}\sum\nolimits_{n=0}^{\infty}\frac{(-1)^nx^{\delta
n}}{\Gamma(\delta n+1)}=-\frac{d}{dx}E_{\delta}(-x^{\delta}),\ \ \
x\ge0,\eqno(2)
$$
where $E_{\delta}(z)$ is the Mittag-Leffler function with index
$\delta$ that is defined as the power series
$$
E_{\delta}(z)=\sum\nolimits_{n=0}^{\infty}\frac{z^n}{\Gamma(\delta
n+1)},\ \ \ \delta>0,\ z\in\mathbb{Z}.
$$
Here $\Gamma(s)$ is Euler's gamma-function,
$$
\Gamma(s)=\int_{0}^{\infty}z^{s-1}e^{-z}dz,\ \ \ s>0.
$$
The distribution function corresponding to density (2) will be
denoted $F_{\delta}^{ML}(x)$.

With $\delta=1$, the Mittag-Leffler distribution turns into the
standard exponential distribution, that is, $F_1^{ML}(x)=
[1-e^{-x}]\mathbf{1}(x\ge 0)$, $x\in\mathbb{R}$ (here and in what
follows the symbol $\mathbf{1}(C)$ denotes the indicator function of
a set $C$). But with $\delta<1$ the Mittag-Leffler distribution
density has the heavy power-type tail: from the well-known
asymptotic properties of the Mittag-Leffler function it can be
deduced that
$$
f_\delta^{ML}(x)\sim \frac{\sin(\delta\pi)\Gamma(\delta+1)}{\pi
x^{\delta+1}}
$$
as $x\to\infty$, see, e. g., \cite{Kilbas2014}.

It is well-known that the Mittag-Leffler distribution is stable with
respect to geometric summation (or {\it geometrically stable}). This
means that if $X_1,X_2,\ldots$ are independent random variables and
$N_p$ is the random variable independent of $X_1,X_2,\ldots$ and
having the geometric distribution
$$
{\sf P}(N_p=n)=p(1-p)^{n-1},\ \ \ n=1,2,\ldots,\ \ \
p\in(0,1),\eqno(3)
$$
then for each $p\in(0,1)$ there exists a constant $a_p>0$ such that
$a_p\big(X_1+\ldots+X_{N_p}\big)\Longrightarrow M_{\delta}$ as $p\to
0$, see, e. g., \cite{Bunge1996} or \cite{KlebanovRachev1996} (the
symbol $\Longrightarrow$ hereinafter denotes convergence in
distribution). Moreover, as far ago as in 1965 it was shown by
I.~Kovalenko \cite{Kovalenko1965} that the distributions with
Laplace transforms (1) are the only possible limit laws for the
distributions of appropriately normalized geometric sums of the form
$a_p\big(X_1+\ldots+X_{N_p}\big)$ as $p\to0$, where $X_1,X_2,\ldots$
are independent identically distributed {\it nonnegative} random
variables and $N_p$ is the random variable with geometric
distribution (3) independent of the sequence $X_1,X_2,\ldots$ for
each $p\in(0,1)$. The proofs of this result were reproduced in
\cite{GnedenkoKovalenko1968, GnedenkoKovalenko1989} and
\cite{GnedenkoKorolev1996}. In these books the class of
distributions with Laplace transforms (1) was not identified as the
class of Mittag-Leffler distributions but was called {\it class}
$\mathcal{K}$ after I.~Kovalenko.

Twenty five years later this limit property of the Mittag-Leffler
distributions was re-discovered by A.~Pillai in \cite{Pillai1989,
Pillai1990} who proposed the term {\it Mittag-Leffler distribution}
for the distribution with Laplace transform (1). Perhaps, since the
works \cite{Kovalenko1965, GnedenkoKovalenko1968,
GnedenkoKovalenko1989} were not easily available to probabilists,
the term {\it class $\mathcal{K}$ distribution} did not take roots
in the literature whereas the term {\it Mittag-Leffler distribution}
became conventional.

The Mittag-Leffler distributions are of serious theoretical interest
in the problems related to thinned (or rarefied) homogeneous flows
of events such as renewal processes or anomalous diffusion or
relaxation phenomena, see \cite{WeronKotulski1996,
GorenfloMainardi2006} and the references therein.

\subsection{The Linnik distributions}

In 1953 Yu. V. Linnik \cite{Linnik1953} introduced the class of
symmetric probability distributions defined by the characteristic
functions
$$
\phi_{\alpha}(t)=\frac{1}{1+|t|^{\alpha}},\ \ \
t\in\mathbb{R},\eqno(4)
$$
where $\alpha\in(0,2]$. Later the distributions of this class were
called {\it Linnik distributions} \cite{Kotz2001} or {\it
$\alpha$-Laplace distributions} \cite{Pillai1985}. In this paper we
will keep to the first term that has become conventional. With
$\alpha=2$, the Linnik distribution becomes the Laplace distribution
corresponding to the density
$$
f^{\Lambda}(x)=\textstyle{\frac12}e^{-|x|},\ \ \
x\in\mathbb{R}.\eqno(5)
$$
A random variable with Laplace density (5) and its distribution
function will be denoted $\Lambda$ and $F^{\Lambda}(x)$,
respectively.

The Linnik distributions possess many interesting analytic
properties such as unimodality \cite{Laha1961} and infinite
divisibility \cite{Devroye1990}, existence of an infinite peak of
the density for $\alpha\le1$ \cite{Devroye1990}, etc. However,
perhaps, most often Linnik distributions are recalled as examples of
geometric stable distributions.

A random variable with the Linnik distribution with parameter
$\alpha$ will be denoted $L_{\alpha}$. Its distribution function and
density will be denoted $F_{\alpha}^{L}$ and $f_{\alpha}^{L}$,
respectively. As this is so, from (4) and (5) it follows that
$F_2^{L}(x)\equiv F^{\Lambda}(x)$, $x\in\mathbb{R}$.

\section{Basic notation and auxiliary results}

Most results presented below actually concern special mixture
representations for probability distributions. However, without any
loss of generality, for the sake of visuality and compactness of
formulations and proofs we will formulate the results in terms of
the corresponding random variables assuming that all the random
variables mentioned in what follows are defined on the same
probability space $(\Omega,\,\mathfrak{A},\,{\sf P})$.

The random variable with the standard normal distribution function
$\Phi(x)$ will be denoted $X$,
$$
{\sf
P}(X<x)=\Phi(x)=\frac{1}{\sqrt{2\pi}}\int_{-\infty}^{x}e^{-z^2/2}dz,\
\ \ \ x\in\mathbb{R}.
$$
Let $\Psi(x)$, $x\in\mathbb{R}$, be the distribution function of the
maximum of the standard Wiener process on the unit interval,
$\Psi(x)=2\Phi\big(\max\{0,x\}\big)-1$, $x\in\mathbb{R}$. It is easy
to see that $\Psi(x)={\sf P}(|X|<x)$. Therefore, sometimes $\Psi(x)$
is said to determine the {\it half-normal} or {\it folded normal}
distribution.

Throughout the paper the symbol $\eqd$ will denote the coincidence
of distributions.

The distribution function and density of the strictly stable
distribution with the characteristic exponent $\alpha$ and shape
parameter $\theta$ defined by the characteristic function
$$
\mathfrak{f}_{\alpha,\theta}(t)=\exp\big\{-|t|^{\alpha}\exp\{-{\textstyle\frac12}i\pi\theta\alpha\mathrm{sign}t\}\big\},\
\ \ \ t\in\r,\eqno(6) 
$$
with $0<\alpha\le2$, $|\theta|\le\min\{1,\frac{2}{\alpha}-1\}$, will
be denoted by $G_{\alpha,\theta}(x)$ and $g_{\alpha,\theta}(x)$,
respectively (see, e. g., \cite{Zolotarev1983}). Any random variable
with the distribution function $G_{\alpha,\theta}(x)$ will be
denoted $S_{\alpha,\theta}$.

From (6) it follows that the characteristic function of a symmetric
($\theta=0$) strictly stable distribution has the form
$$
\mathfrak{f}_{\alpha,0}(t)=e^{-|t|^{\alpha}},\ \ \ t\in\r. \eqno(7) 
$$
From (7) it is easy to see that $S_{2,0}\eqd\sqrt{2}X$.

\smallskip

{\sc Lemma 1}. {\it Let $\alpha\in(0,2]$, $\alpha'\in(0,1]$. Then}
$$
S_{\alpha\alpha',0}\eqd S_{\alpha,0}S_{\alpha',1}^{1/\alpha}
$$
{\it where the random variables on the right-hand side are
independent.}

\smallskip

{\sc Proof}. See, e. g., theorem 3.3.1 in \cite{Zolotarev1983}.

\smallskip

{\sc Corollary 1.} {\it A symmetric strictly stable distribution
with the characteristic exponent $\alpha$ is a scale mixture of
normal laws in which the mixing distribution is the one-sided
strictly stable law $(\theta=1)$ with the characteristic exponent
$\alpha/2$}:
$$
G_{\alpha,0}(x)=\int_{0}^{\infty}\Phi\big(x/\sqrt{2z}\big)dG_{\alpha/2,1}(z),\
\ \ x\in\r.
$$

\smallskip

In terms of random variables the statement of corollary 1 can be
written as
$$
S_{\alpha,0}\eqd X\sqrt{2S_{\alpha/2,1}} \eqno(8)
$$
with the random variables on the right-hand side being independent.

\smallskip

Let $\gamma>0$. The distribution of the random variable
$W_{\gamma}$:
$$
{\sf
P}\big(W_{\gamma}<x\big)=\big[1-e^{-x^{\gamma}}\big]\mathbf{1}(x\ge
0),\ \ \ x\in\mathbb{R},
$$
is called the {\it Weibull distribution} with shape parameter
$\gamma$. It is obvious that $W_1$ is the random variable with the
standard exponential distribution: ${\sf
P}(W_1<x)=\big[1-e^{-x}\big]{\bf 1}(x\ge0)$. The Weibull
distribution with $\gamma=2$, that is, ${\sf
P}(W_2<x)=\big[1-e^{-x^2}\big]{\bf 1}(x\ge0)$ is called the Rayleigh
distribution.

It is easy to see that if $\gamma>0$ and $\gamma'>0$, then ${\sf
P}(W_{\gamma'}^{1/\gamma}\ge x)={\sf P}(W_{\gamma'}\ge
x^{\gamma})=e^{-x^{\gamma\gamma'}}={\sf P}(W_{\gamma\gamma'}\ge x)$,
$x\ge 0$, that is, for any $\gamma>0$ and $\gamma'>0$
$$
W_{\gamma\gamma'}\eqd W_{\gamma'}^{1/\gamma}.\eqno(9) 
$$
In particular, for any $\gamma>0$ we have $W_{\gamma}\eqd
W_1^{1/\gamma}$.

It can be shown that each Weibull distribution with parameter
$\gamma\in(0,1]$ is a mixed exponential distribution. In order to
prove this we first make sure that each Weibull distribution with
parameter $\gamma\in(0,2]$ is a scale mixture of the Rayleigh
distributions.

For $\alpha\in(0,1]$ denote $V_{\alpha}=S_{\alpha,1}^{-1}$, where
$S_{\alpha,1}$ is a random variable with one-sided strictly stable
density $g_{\alpha,1}(x)$.

\smallskip

{\sc Lemma 2.} {\it For any $\gamma\in(0,2]$ we have}
$$
W_{\gamma}\eqd W_2 \sqrt{V_{\gamma/2}},\eqno(10) 
$$
{\it where the random variables on the right-hand side of $(10)$ are
independent.}

\smallskip

{\sc Proof}. Write relation (8) with $\alpha$ replaced by $\gamma$
in terms of characteristic functions with the account of (7):
$$
e^{-|t|^{\gamma}}=\int_{0}^{\infty}\exp\{-t^2z\}g_{\gamma/2,1}(z)dz,
\ \ \ t\in\mathbb{R}.\eqno(11) %
$$
Formally letting $|t|=x$ in (11), where $x\ge0$ is an arbitrary
nonnegative number, we obtain
$$
{\sf
P}(W_{\gamma}>x)=e^{-x^{\gamma}}=\int_{0}^{\infty}\exp\{-x^2z\}g_{\gamma/2,1}(z)dz.\eqno(12) 
$$
At the same time it is obvious that if $W_2$ and $S_{\gamma/2,1}$
are independent, then
$$
{\sf P}\big(W_2 \sqrt{V_{\gamma/2}}>x\big)={\sf
P}\big(W_2>x\sqrt{S_{\gamma/2,1}}\big)=
\int_{0}^{\infty}\exp\{-x^2z\}g_{\gamma/2,1}(z)dz.\eqno(13) 
$$
Since the right-hand sides of (12) and (13) coincide identically in
$x\ge0$, the left-hand sides of these relations coincide as well.
The lemma is proved.

\smallskip

{\sc Lemma 3}. {\it For any $\gamma\in(0,1]$, the Weibull
distribution with parameter $\gamma$ is a mixed exponential
distribution}:
$$
W_{\gamma}\eqd W_1 V_{\gamma}.\eqno(14) 
$$
{\it where the random variables on the right-hand side of $(14)$ are
independent.}

\smallskip

{\sc Proof.} It is easy to see that ${\sf P}(W_1^{1/\gamma}\ge
x)={\sf P}(W_1\ge x^{\gamma})=e^{-x^{\gamma}}={\sf P}(W_{\gamma}\ge
x)$, $x\ge 0$, that is,
$$
W_{\gamma}\eqd W_1^{1/\gamma}.\eqno(15) 
$$
From $(15)$ it follows that $W_2\eqd\sqrt{W_1}$. Therefore, from
lemma 2 it follows that for $\gamma\in(0,2]$ we have
$$
W_{\gamma}\eqd W_2 \sqrt{V_{\gamma/2}}\eqd\sqrt{W_1 V_{\gamma/2}}
$$
or, with the account of $(15)$,
$$
W_{\gamma/2}\eqd W_{\gamma}^2\eqd W_1 V_{\gamma/2}.
$$
Re-denoting $\gamma/2\longmapsto \gamma\in(0,1]$, we obtain the
desired assertion.

\smallskip

In \cite{Devroye1990} the following statement was proved. Here its
formulation is extended with the account of (9).

\smallskip

{\sc Lemma 4} \cite{Devroye1990}. {\it For any $\alpha\in(0,2]$, the
Linnik distribution with parameter $\alpha$ is a scale mixture of a
symmetric stable distribution with the Weibull mixing distribution
with parameter $\alpha/2$, that is,}
$$
L_{\alpha}\eqd S_{\alpha,0}W_{\alpha}\eqd
S_{\alpha,0}\sqrt{W_{\alpha/2}},
$$
{\it where the random variables on the right-hand side are
independent}.

\smallskip

{\sc Lemma 5.} {\it For any $\delta\in(0,1]$, the Mittag-Leffler
distribution with parameter $\delta$ is a scale mixture of a
one-sided stable distribution with the Weibull mixing distribution
with parameter $\delta/2$, that is,}
$$
M_{\delta}\eqd S_{\delta,1}W_{\delta}\eqd
S_{\delta,1}\sqrt{W_{\delta/2}},
$$
{\it where the random variables on the right-hand side are
independent}.

\smallskip

{\sc Proof}. This statement has already become folklore. For the
purpose of convenience we give its elementary proof without any
claims for priority. Let $S_{\delta,1}$ be a positive strictly
stable random variable. As is known, its Laplace transform is
$\psi(s)={\sf E}e^{-sS_{\delta,1}}=e^{-s^{\delta}}$, $s\ge0$. Then
with the account of (9) by the Fubini theorem the Laplace transform
of the product $S_{\delta,1}W_{\delta}$ is
$$
{\sf E}\exp\{-sS_{\delta,1}W_{\delta}\}= {\sf
E}\exp\{-sS_{\delta,1}W_1^{1/\delta}\}={\sf E}{\sf
E}\big(\exp\{-sS_{\delta,1}W_1^{1/\delta}\}\big|W_1\big)=
\int_{0}^{\infty}e^{-(sz^{1/\delta})^{\delta}}e^{-z}dz=
$$
$$
=\int_{0}^{\infty}e^{-z(s^{\delta}+1)}dz=\frac{1}{1+s^{\delta}}={\sf
E}e^{-sM_{\delta}},\ \ \ s\ge0.
$$
The lemma is proved.

\smallskip

Let $\rho\in(0,1)$. In \cite{Kozubowski1998} it was demonstrated
that the function
$$
f_{\rho}^{K}(x)=\frac{\sin(\pi\rho)}{\pi\rho[x^2+2x\cos(\pi\rho)+1]},\
\ \ x\in(0,\infty),\eqno(16)
$$
is a probability density on $(0,\infty)$. Let $K_{\rho}$ be a random
variable with density (16).

\smallskip

{\sc Lemma 6} \cite{Kozubowski1998}. {\it Let $0<\delta<\delta'\le1$
and $\rho=\delta/\delta'<1$. Then}
$$
M_{\delta}\eqd M_{\delta'}K_{\rho}^{1/\delta}
$$
{\it where the random variables on the right-hand side are
independent.}

\smallskip

With $\delta'=1$ we have

\smallskip

{\sc Corollary 2} \cite{Kozubowski1998}. {\it Let $0<\delta<1$. Then
the Mittag-Leffler distribution with parameter $\delta$ is mixed
exponential}:
$$
M_{\delta}\eqd K_{\delta}^{1/\delta}W_1
$$
{\it where the random variables on the right-hand side are
independent.}

\smallskip

Let $0<\alpha<\alpha'\le2$. In \cite{KotzOstrovskii1996} it was
shown that the function
$$
f_{\alpha,\alpha'}^{Q}(x)= \frac{\alpha'\sin(\pi\alpha/\alpha')
x^{\alpha-1}}{\pi[1+x^{2\alpha}+2x^{\alpha}\cos(\pi\alpha/\alpha')]},\
\ \ x>0,\eqno(17)
$$
is a probability density on $(0,\infty)$. Let $Q_{\alpha,\alpha'}$
be a random variable whose probability density is
$f_{\alpha,\alpha'}^{Q}(x)$.

\smallskip

{\sc Lemma 7} \cite{KotzOstrovskii1996}. {\it Let
$0<\alpha<\alpha'\le2$. Then}
$$
L_{\alpha}\eqd L_{\alpha'}Q_{\alpha,\alpha'},
$$
{\it where the random variables on the right-hand side are
independent}.

\smallskip

With $\alpha'=2$ we have

\smallskip

{\sc Corollary 3} \cite{KotzOstrovskii1996}. {\it Let $0<\alpha<2$.
Then the Linnik distribution with parameter $\alpha$ is a scale
mixture of Laplace distributions corresponding to density $(5)$}:
$$
L_{\alpha}\eqd \Lambda Q_{\alpha,2}
$$
{\it where the random variables on the right-hand side are
independent.}

\smallskip

For the sake of completeness, we will demonstrate that the Weibull
distributions possess the same property as the Linnik and
Mittag-Leffler distributions presented in lemmas 6 and 7: any
distribution of the corresponding class can be represented as a
scale mixture of a distribution from the same class with larger
parameter.

Relation (14) implies the following statement generalizing lemmas 2
and 3 and stating that the Weibull distribution with an arbitrary
positive shape parameter $\gamma$ is a scale mixture of a Weibull
distribution with an arbitrary positive shape parameter
$\gamma'>\gamma$.

\smallskip

{\sc Lemma 8}. {\it Let $\gamma'>\gamma>0$ be arbitrary numbers.
Then
$$
W_{\gamma}\eqd W_{\gamma'}\cdot V_{\alpha}^{1/\gamma'},
$$
where $\alpha=\gamma/\gamma'\in(0,1)$ and the random variables on
the right-hand side are independent}.

\smallskip

{\sc Proof}. In lemma 3 we showed that a Weibull distribution with
parameter $\alpha\in(0,1]$ is a mixed exponential distribution.
Indeed, from (14) it follows that
$$
e^{-x^{\alpha}}={\sf P}(W_{\alpha}>x)={\sf P}(W_1>S_{\alpha,1}x)=
\int_{0}^{\infty}e^{-zx}g_{\alpha,1}(z)dz,\ \ \ x\ge0.
$$
Therefore, for any $\gamma'>\gamma>0$, denoting
$\alpha=\gamma/\gamma'$ (as this is so, $\alpha\in(0,1)$), for any
$x\in\mathbb{R}$ we obtain
$$
{\sf P}(W_{\gamma}>x)=e^{-x^{\gamma}}=e^{-x^{\gamma'\alpha}}={\sf
P}(W_{\alpha}>x^{\gamma'})={\sf P}(W_1>S_{\alpha,1}x^{\gamma'})=
$$
$$
=\int_{0}^{\infty}e^{-zx^{\gamma'}}g_{\alpha,1}(z)dz=\int_{0}^{\infty}{\sf
P}\big(W_{\gamma'}>xz^{1/\gamma'}\big)g_{\alpha,1}(z)dz={\sf
P}(W_{\gamma'}\cdot V_{\alpha}^{1/\gamma'}>x),
$$
The lemma is proved.

\smallskip

It should be noted that if $0<\gamma<\gamma'<2$, then the assertion
of lemma 8 directly follows from theorem 3.3.1 of
\cite{Zolotarev1983} due to the formal coincidence of the
characteristic function of a strictly stable law and the
complementary Weibull distribution function (see the proof of lemma
2).

\section{Representation of the Linnik distribution as a scale
mixture of normal or Laplace distributions and related results}

In all the products of random variables mentioned below the
multipliers are assumed independent.

The following statement presents the main result of this paper.

\smallskip

{\sc Theorem 1.} {\it Let $\alpha\in(0,2]$, $\alpha'\in(0,1]$. Then}
$$
L_{\alpha\alpha'}\eqd S_{\alpha,0}M_{\alpha'}^{1/\alpha}.
$$

\smallskip

{\sc Proof}. From lemma 4 we have
$$
L_{\alpha\alpha'}\eqd
S_{\alpha\alpha',0}\sqrt{W_{\alpha\alpha'/2}}.\eqno(18)
$$
Continuing (18) with the account of lemma 1, we obtain
$$
L_{\alpha\alpha'}\eqd
S_{\alpha,0}S_{\alpha',1}^{1/\alpha}\sqrt{W_{\alpha\alpha'/2}}.\eqno(19)
$$
From (9) and lemma 5 it follows that
$$
S_{\alpha',1}^{1/\alpha}\sqrt{W_{\alpha\alpha'/2}}\eqd
S_{\alpha',1}^{1/\alpha}W_{\alpha\alpha'}\eqd
(S_{\alpha',1}W_{\alpha'})^{1/\alpha}\eqd M_{\alpha'}^{1/\alpha}.
$$
The theorem is proved.

\smallskip

As far as we know, the following result has never been explicitly
presented in the literature in full detail although the property of
the Linnik distribution to be a normal scale mixture is something
almost obvious.

\smallskip

{\sc Corollary 4}. {\it For each $\alpha\in(0,2]$, the Linnik
distribution with parameter $\alpha$ is the scale mixture of
zero-mean normal laws with mixing Mittag-Leffler distribution with
twice less parameter $\alpha/2$}:
$$
L_{\alpha}\eqd X\sqrt{2M_{\alpha/2}},
$$
{\it where the random variables on the right-hand side are
independent}.

\smallskip

It should be noted that in this case representation (19) takes the
form
$$
L_{\alpha}\eqd X\sqrt{2S_{\alpha/2,1}W_{\alpha/2}}.\eqno(20)
$$
From lemma 3 we have
$$
W_{\alpha/2}\eqd \frac{W_1}{S'_{\alpha/2,1}}.
$$
Hence, from (20) it follows that
$$
L_{\alpha}\eqd X\sqrt{2W_1\frac{S_{\alpha/2,1}}{S'_{\alpha/2,1}}}
$$
where the independent random variables $S_{\alpha/2,1}$ and
$S'_{\alpha/2,1}$ have one and the same one-sided strictly stable
distribution with characteristic exponent $\alpha/2$ and are
independent of the exponentially distributed random variable $W_1$.
It is well known that
$$
X\sqrt{2W_1}\eqd \Lambda\eqno(21)
$$
(see, e. g., the example on p. 272 of \cite{BeningKorolev2002}).
Therefore we obtain one more mixture representation for the Linnik
distribution.

\smallskip

{\sc Theorem 2}. {\it For each $\alpha\in(0,2]$, the Linnik
distribution with parameter $\alpha$ is the scale mixture of the
Laplace laws corresponding to density $(5)$ with mixing distribution
being that of the ratio of two independent random variables having
one and the same one-sided strictly stable distribution with
characteristic exponent $\alpha/2$}:
$$
L_{\alpha}\eqd \Lambda\sqrt{\frac{S_{\alpha/2,1}}{S'_{\alpha/2,1}}},
$$
{\it where the random variables on the right-hand side are
independent}.

\smallskip

It is easy to see that scale mixtures of Laplace distribution (5)
are identifiable, that is, if
$$
\Lambda Y\eqd\Lambda Y'
$$
where $Y$ and $Y$ are nonnegative random variables independent of
$\Lambda$, then $Y\eqd Y'$. Indeed, with the account of (21), the
last relation turns into
$$
X\sqrt{2W_1 Y^2}\eqd X\sqrt{2W_1 (Y')^2},\eqno(22)
$$
where the random the multipliers on both sides are independent. But,
as is known, scale mixtures of zero-mean normals are identifiable
(see \cite{Teicher1961}). Therefore, (22) implies that
$$
W_1 Y^2\eqd W_1(Y')^2.\eqno(23)
$$
The complementary mixed exponential distribution functions of the
random variables related by (23) are the Laplace transforms of $Y^2$
and $(Y')^2$, respectively. Relation (23) means that these Laplace
transforms identically coincide:
$$
\int_{0}^{\infty}e^{-sz}d{\sf
P}(Y^2<z)\equiv\int_{0}^{\infty}e^{-sz}d{\sf P}\big((Y')^2<z\big),\
\ \ s\ge0.
$$
Hence, the distributions of the random variables $Y^2$ and $(Y')^2$
coincide and hence, the distributions of $Y$ and $Y'$ coincide as
well since these random variables were originally assumed
nonnegative.

Comparing the statement of theorem 2 with the assertion of corollary
3 with the account of identifiability of scale mixtures of Laplace
distributions (5) we arrive at the relation
$$
Q_{\alpha,2}\eqd\sqrt{\frac{S_{\alpha/2,1}}{S'_{\alpha/2,1}}}.\eqno(24)
$$
The combination of (17) and (24) gives one more, possibly simpler,
proof of the following by-product result concerning the properties
of stable distributions obtained in \cite{Dovidio2010}. This result
offers an explicit representation for the density of the ratio of
two independent stable random variables in terms of elementary
functions although with the exception of one case, the L{\'e}vy
distribution $(\alpha=\frac12)$, such representations for the
densities of nonnegative stable random variables themselves do not
exist.

\smallskip

{\sc Corollary 5}. {\it Let $S_{\alpha,1}$ and $S'_{\alpha,1}$ be
two independent random variables having one and the same one-sided
strictly stable distribution with characteristic exponent
$\alpha\in(0,1)$. Then $S_{\alpha,1}/S'_{\alpha,1}\eqd
K_{\alpha}^{1/\alpha}\eqd Q^2_{2\alpha,2}$, that is, the probability
density $p_{\alpha}(x)$ of the ratio $S_{\alpha,1}/S'_{\alpha,1}$
has the form}
$$
p_{\alpha}(x)=f_{\alpha,1}^{Q}(x)=\frac{\sin(\pi\alpha)x^{\alpha-1}}{\pi[1+x^{2\alpha}+2x^{\alpha}\cos(\pi\alpha)]},\
\ \ x>0.
$$

\smallskip

The {\sc proof} immediately follows from the observation that
$p_{\alpha}(x)=(2\sqrt{x})^{-1}f_{2\alpha,2}^{Q}(\sqrt{x})=f_{\alpha,1}^{Q}(x)$.

\smallskip

As concerns the Mittag-Leffler distribution, from lemmas 3 and 5 we
obtain the following statement analogous to theorem 2.

\smallskip

{\sc Theorem 3}. {\it For each $\delta\in(0,1]$, the Mittag-Leffler
distribution with parameter $\delta$ is the mixed exponential
distribution with mixing distribution being that of the ratio of two
independent random variables having one and the same one-sided
strictly stable distribution with characteristic exponent $\delta$}:
$$
M_\delta\eqd W_1\frac{S_{\delta,1}}{S'_{\delta,1}},
$$
{\it where the random variables on the right-hand side are
independent}.

\smallskip

From theorem 3 and corollary 5 we obtain the following
representation of the Mittag-Leffler distribution function
$F_{\delta}^{ML}(x)$:
$$
F_{\delta}^{ML}(x)=1-\frac{\sin(\pi\delta)}{\pi}\int_{0}^{\infty}\frac{
z^{\delta-1}e^{-zx}dz}{1+z^{2\delta}+2z^{\delta}\cos(\pi\delta)},\ \
\ x>0.\eqno(25)
$$

Representation for the Linnik distribution similar to (25) was
obtained in \cite{Kozubowski1998}.

Using lemmas 1, 4 and 5 it  it is possible to obtain more product
representations for the Mittag-Leffler- and Linnik-distributed
random variables and hence, more mixture representations for these
distributions.

\section{Convergence of the distributions of random sums to the
Linnik distribution}

Product representations for the random variables with the Linnik and
Mittag-Leffler distributions obtained in the previous works were
aimed at the construction of convenient algorithms for the computer
generation of pseudo-random variables with these distributions. The
mixture representation for the Linnik distribution as a scale
mixture of normals obtained in corollary 4 opens the way for the
construction in this section of a random-sum central limit theorem
with the Linnik distribution as the limit law. Moreover, in this
version of the random-sum central limit theorem the Mittag-Leffler
distribution {\it must} be the limit law for the normalized number
of summands.

Recall that the symbol $\Longrightarrow$ denotes the convergence in
distribution.

Consider a sequence of independent identically distributed random
variables $X_1,X_2,\ldots$ defined on the probability space
$(\Omega,\, \mathfrak{A},\,{\sf P})$. Assume that ${\sf E}X_1=0$,
$0<\sigma^2={\sf D}X_1<\infty$. For $n\in\mathbb{N}$ denote
$S^*_n=X_1+\ldots+X_n$. Let $Z_1,Z_2,\ldots$ be a sequence of
integer-valued nonnegative random variables defined on the same
probability space so that for each $n\ge1$ the random variable $Z_n$
is independent of the sequence $X_1,X_2,\ldots$ For definiteness, in
what follows we assume that $\sum\nolimits_{j=1}^0=0$.

Recall that a random sequence $Z_1,Z_2,\ldots$ is said to infinitely
increase in probability ($Z_n\pto\infty$), if ${\sf P}(Z_n\le
m)\longrightarrow 0$ as $n\to\infty$ for any $m\in(0,\infty)$.

The proof of the main result of this section is based on the
following version of the random-sum central limit theorem.

\smallskip

{\sc Lemma 10.} {\it Assume that the random variables
$X_1,X_2,\ldots$ and $Z_1,Z_2,\ldots$ satisfy the conditions
specified above and, moreover, let $Z_n\pto\infty$ as $n\to\infty$.
A distribution function $F(x)$ such that
$$
{\sf P}\Big(\frac{S^*_{Z_n}}{\sigma\sqrt{n}}<x\Big) \Longrightarrow
F(x)
$$
as $n\to\infty$ exists if and only if there exists a distribution
function $H(x)$ satisfying the conditions
$$
H(0)=0,\ \ \
F(x)=\int_{0}^{\infty}\Phi\Big(\frac{x}{\sqrt{y}}\Big)dH(y),\ \
x\in\mathbb{R},
$$
and ${\sf P}(Z_n<nx)\Longrightarrow H(x)$ $(n\to\infty)$. }

\smallskip

{\sc Proof}. This statement is a particular case of a result proved
in \cite{Korolev1994}, also see theorem 3.3.2 in
\cite{GnedenkoKorolev1996}.

\smallskip

The following theorem gives a criterion (that is, {\it necessary and
sufficient} conditions) of the convergence of the distributions of
random sums of independent identically distributed random variables
with {\it finite} variances to the Linnik distribution.

\smallskip

{\sc Theorem 4.} {\it Let $\alpha\in(0,2]$. Assume that the random
variables $X_1,X_2,\ldots$ and $Z_1,Z_2,\ldots$ satisfy the
conditions specified above and, moreover, let $Z_n\pto\infty$ as
$n\to\infty$. Then the distributions of the normalized random sums
$S^*_{Z_n}$ converge to the Linnik law with parameter $\alpha$, that
is,
$$
{\sf P}\Big(\frac{S^*_{Z_n}}{\sigma\sqrt{n}}<x\Big) \Longrightarrow
F_{\alpha}^{L}(x)
$$
as $n\to\infty$, if and only if
$$
\frac{Z_n}{n}\Longrightarrow 2M_{\alpha/2}\ \ \ (n\to\infty).
$$
}

\smallskip

{\sc Proof}. This statement is a direct consequence of corollary 4
and lemma 10 with $H(x)=F_{\alpha/2}^{ML}(x/2)$.

\smallskip

The convergence of the distributions of the normalized indices
$Z_n/n$ to the Mittag-Leffler distribution $F_{\delta}^{ML}$ is the
main condition in theorem 4. Now we will give two examples of the
situation where this condition can hold. The first example is
trivial and is based on the geometric stability of the
Mittag-Leffler distribution. The second example relies on a useful
general construction of nonnegative integer-valued random variables
which, under an appropriate normalization, converge to a given
nonnegative (not necessarily discrete) random variable, whatever the
latter is.

\smallskip

{\sc Example 1}. Let $\delta\in(0,1)$ be arbitrary. For every
$n\in\mathbb{N}$ let $N_{1/n}$ be a random variable having the
geometric distribution (3) with $p=\frac1n$ independent of the
sequence $Y_1,Y_2,\ldots$ of independent identically distributed
nonnegative random variables such that
$$
n^{-1/\delta}\sum\nolimits_{j=1}^{N_{1/n}}Y_j\Longrightarrow
2M_{\delta}\eqno(26)
$$
as $n\to\infty$. To provide (26), the distributions of the random
variables $Y_1,Y_2,\ldots$ should belong to the domain of the normal
attraction of the one-sided strictly stable law with characteristic
exponent $\delta$. As $Z_n$ for each $n\in\mathbb{N}$ take
$$
Z_n=\bigg[n^{1-1/\delta}\sum\nolimits_{j=1}^{N_{1/n}}Y_j\bigg],
$$
where square brackets denote the integer part. Then
$$
\frac{Z_n}{n}=n^{-1/\delta}\sum\nolimits_{j=1}^{N_{1/n}}Y_j-\frac{1}{n}\bigg\{n^{1-1/\delta}\sum\nolimits_{j=1}^{N_{1/n}}Y_j\bigg\},\eqno(27)
$$
where curly braces denote the fractional part. Since the second term
on the right-hand side of (27) obviously tends to zero in
probability, from (26) it follows that $Z_n/n\Longrightarrow
2M_{\delta}$ as $n\to\infty$.

\smallskip

{\sc Example 2.} In the book \cite{GnedenkoKorolev1996} it was
proposed to model the evolution of non-homogeneous chaotic
stochastic processes, in particular, the dynamics of financial
markets by compound doubly stochastic Poisson processes (compound
Cox processes). This approach got further grounds and development,
say, in \cite{BeningKorolev2002, Zeifman2015}. According to this
approach the flow of informative events, each of which generates the
next observation, is described by the stochastic point process
$P(Y(t))$ where $P(t)$, $t\geq0$, is a homogeneous Poisson process
with unit intensity and $Y(t)$, $t\geq0$, is a random process
independent of $P(t)$ and possessing the properties: $Y(0)=0$, ${\sf
P}(Y(t)<\infty)=1$ for any $t>0$, the trajectories $Y(t)$ are
non-decreasing and right-continuous. The process $P(Y(t))$,
$t\geq0$, is called a doubly stochastic Poisson process (Cox
process) \cite{Grandell1976}.

Within this model, for each $t$ the distribution of the random
variable $P(Y(t))$ is mixed Poisson. For vividness, consider the
case where in the model under consideration the parameter $t$ is
discrete: $Y(t)=Y(n)=Y_n$, $n\in\mathbb{N}$, where $\{Y_n\}_{n\ge1}$
is an infinitely increasing sequence of nonnegative random variables
such that $Y_{n+1}(\omega)\ge Y_{n}(\omega)$ for any
$\omega\in\Omega$, $n\ge1$. Here the asymptotics $n\to\infty$ may be
interpreted as that the intensity of the flow of informative events
is (infinitely) large.

Assume that the random variable $M_{\delta}$ is independent of the
standard Poisson process $P(t)$, $t\ge0$. For each natural $n$ take
$Y_n=2nM_{\delta}$. Respectively, let $Z_n=P(Y_n)=P(2nM_{\delta})$,
$n\ge1$. It is obvious that the random variable $Z_n$ so defined has
the mixed Poisson distribution
$$
{\sf P}(Z_n=k)={\sf
P}\big(P(2nM_{\delta})=k\big)=\int_{0}^{\infty}e^{-2nz}\frac{(2nz)^k}{k!}dF_{\delta}^{ML}(z)\
\ \ k=0,1,...
$$
This random variable $Z_n$ can be interpreted as the number of
events registered up to time $n$ in the Poisson process with the
stochastic intensity distributed as $2M_{\delta}$.

Denote $A_n(z)={\sf P}(Z_n<2nz)$, $z\ge0$ ($A_n(z)=0$ for $z<0$). It
is easy to see that $A_n(z)\Longrightarrow F_{\delta}^{ML}(z)$ as
$n\to\infty$. Indeed, as is known, if $\Pi(x;\ell)$ is the Poisson
distribution function with the parameter $\ell>0$ and $E(x;c)$ is
the distribution function with a single unit jump at the point
$c\in\r$, then $\Pi(\ell x;\ell)\Longrightarrow E(x;1)$ as
$\ell\to\infty$. Since for $x\in\r$
$$
A_n(x)=\int_{0}^{\infty}\Pi(2n x; 2n z)dF_{\delta}^{ML}(z),
$$
then by the Lebesgue dominated convergence theorem, as $n\to\infty$,
we have
$$
A_n(x)\Longrightarrow\int_{0}^{\infty}E(x/z;1)dF_{\delta}^{ML}(z)=
\int_{0}^{x}dF_{\delta}^{ML}(z)=F_{\delta}^{ML}(x),
$$
that is, the random variables $Z_n$ defined above satisfy the
condition of theorem 4. Moreover, $Z_n\pto\infty$ as $n\to\infty$
since ${\sf P}(M_{\delta}=0)=0$.

\smallskip

The authors thank Dmitri Alekritskii who carefully read the draft
text of the paper and pointed out some inaccuracies.

\renewcommand{\refname}{References}


\begin{thebibliography}{99}

\bibitem{BeningKorolev2002}
{\it Bening V., Korolev V.} Generalized Poisson Models and their
Applications in Insurance and Finance. -- Utrecht: VSP, 2002.

\bibitem{Bunge1996} {\it Bunge J.} Compositions semigroups and random stability // Annals of
Probability, 1996. Vol. 24. P. 1476--1489.

\bibitem{Devroye1990} {\it Devroye, L.} A note on Linnik's
distribution // Statistics and Probability Letters, 1990. Vol. 9. P.
305--306.

\bibitem{Dovidio2010} {\it D'ovidio, M.} Explicit solutions to fractional
diffusion equations via generalized gamma convolution // Electronic
Communications in Probability, 2010. Vol. 15. P. 457--474.

\bibitem{ErdoganOstrovskii1997} {\it Erdo\v{g}an M. B., Ostrovskii I. V.}
Analytic and asymptotic properties of generalized Linnik probability
densities // Journal of Mathematical Analysis and Applications,
1998. Vol. 217. P. 555--578.

\bibitem{ErdoganOstrovskii1998} {\it Erdo\v{g}an M. B., Ostrovskii I. V.}
On mixture representation of the Linnik density // Journal of the
Australian Mathematical Society. Ser. A, 1998. Vol. 64. P. 317--326.

\bibitem{GnedenkoKovalenko1968} {\it Gnedenko B. V., Kovalenko I. N.} Introduction to Queueing
Theory. -- Jerusalem: Israel Program for Scientific Translations,
1968.

\bibitem{GnedenkoKovalenko1989} {\it Gnedenko B. V., Kovalenko I. N.} Introduction to Queueing
Theory. 2nd Edition. -- Boston: Birkhauser, 1989.

\bibitem{GnedenkoKorolev1996} {\it Gnedenko B. V., Korolev V. Yu.} Random Summation:
Limit Theorems and Applications. -- Boca Raton: CRC Press, 1996.

\bibitem{GorenfloMainardi2006} {\it Gorenflo R., Mainardi F.} Continuous time random walk, Mittag-Leffler
waiting time and fractional diffusion: mathematical aspects / Chap.
4 in {\it Klages R., Radons G. and Sokolov I. M.$($Editors$)$.}
Anomalous Transport: Foundations and Applications. -- Weinheim,
Germany: Wiley-VCH, 2008, p. 93--127. Available at:
http://arxiv.org/abs/0705.0797.

\bibitem{Kilbas2014} {\it Gorenflo R., Kilbas A. A., Mainardi F., Rogosin S. V.} Mittag-Leffler Functions,
Related Topics and Applications. -- Berlin-New York: Springer, 2014.

\bibitem{Grandell1976} {\it J. Grandell}. Doubly Stochastic Poisson
Processes. Lecture Notes Mathematics, Vol. 529. --
Berlin-Heidelberg-New York: Springer, 1976.

\bibitem{KlebanovRachev1996} {\it Klebanov L. B., Rachev S. T.} Sums of a random number of random
variables and their approximations with $\varepsilon$-accompanying
infinitely divisible laws // Serdica, 1996. Vol. 22. P. 471--498.

\bibitem{Korolev1994} {\it Korolev V. Yu.} Convergence of random
sequences with independent random indexes. I // Theory of
Probability and iits Applications, 1994. Vol. 39. No. 2. P.
313--333.

\bibitem{Zeifman2015} {\it Korolev V. Yu, Chertok A. V.,
Korchagin A. Yu, Zeifman A. I.} Modeling high-frequency order flow
imbalance by functional limit theorems for two-sided risk processes
// Applied Mathematics and Computation, 2015. Vol. 253. P. 224--241.

\bibitem{KotzOstrovskii1996} {\it Kotz S., Ostrovskii I. V.} A mixture
representation of the Linnik distribution // Statistics and
Probability Letters, 1996. Vol. 26. P. 61--64.

\bibitem{Kotz2001} {\it Kotz S., Kozubowski T. J., Podgorski K.} The
Laplace Distribution and Generalizations: A Revisit with
Applications to Communications, Economics, Engineering, and Finance.
-- Boston: Birkhauser, 2001.

\bibitem{Kovalenko1965} {\it Kovalenko I. N.} On the class of limit
distributions for rarefied flows of homogeneous events // Litovskii
Matematicheskii Sbornik (Lithuanian Mathematical Journal), 1965,
Vol. 5. No. 4. P. 569--573.

\bibitem{Kozubowski1998} {\it Kozubowski T. J.} Mixture representation
of Linnik distribution revisited // Statistics and Probability
Letters, 1998. Vol. 38. P. 157--160.

\bibitem{Kozubowski1999} {\it Kozubowski T. J.} Exponential mixture
representation of geometric stable distributions // Annals of the
Institute of Statistical Mathematics, 1999. Vol. 52. No. 2. P.
231--238.

\bibitem{Laha1961} {\it Laha R. G.} On a class of unimodal
distributions // Proceedings of the American Mathematical Society,
1961. Vol. 12. P. 181--184.

\bibitem{Linnik1953} {\it Linnik Yu. V.} Linear forms and statistical criteria, I,
II // Selected Translations in Mathematical Statistics and
Probability, 1963. Vol. 3. P. 41--90 (Original paper appeared in:
Ukrainskii Matematicheskii Zhournal, 1953. Vol. 5. P. 207--243,
247--290).

\bibitem{Pakes1992} {\it Pakes A. G.} A characterization of gamma
mixtures of stable laws motivated by limit theorems // Statistica
Neerlandica, 1992. Vol. 2-3. P. 209--218.

\bibitem{Pakes1997} {\it Pakes A. G.} Mixture representations for symmetric generalized
Linnik laws // Statistics and Probability Letters, 1998. Vol. 37. P.
213--221.

\bibitem{Pillai1985} {\it Pillai R. N.} Semi-$\alpha$-Laplace
distributions // Communications in Statistical Theory and Methods,
1985. Vol. 14. P. 991-1000.

\bibitem{Pillai1989} {\it Pillai R. N.} Harmonic mixtures and geometric infinite divisibility // Journal of
Indian Statistical Association, 1990. Vol. 28. P. 87--98.

\bibitem{Pillai1990} {\it Pillai R. N.} On Mittag-Leffler functions and related distributions //
Annals of the Institute of Statistical Mathematics, 1990. Vol. 42.
P. 157--161.

\bibitem{Teicher1961} {\it Teicher H.} Identifiability of mixtures // Annals of Mathematical Statistics,
1961. Vol. 32. P. 244--248.

\bibitem{WeronKotulski1996} {\it Weron K., Kotulski M.} On the
Cole-Cole relaxation function and related Mittag-Leffler
distributions // Physica A, 1996. Vol. 232. P. 180--188.

\bibitem{Zolotarev1983} {\it Zolotarev V. M.} One-Dimensional Stable
Distributions. Translation of Mathematical Monographs, Vol. 65. --
Providence, RI: American Mathematical Society, 1986.

\end{thebibliography}
\end{document}